\documentclass[12pt]{elsart}
\usepackage{amsfonts,amsmath}
\usepackage{graphicx}
\usepackage{epic}

\newtheorem{df}[thm]{\bf Definition}
\journal{
}
\date{}
\begin{document}
\begin{frontmatter}
\title
{ Irreducibility criterion for the set of two matrices}
%\thanksref{dec}

\thanks[dec]{The  author would like to thank the Max-Planck-Institute
of Mathematics, Bonn and the Institute of Applied Mathematics,
University of Bonn for the hospitality. The partial financial
support by the DFG project 436 UKR 113/87 is gratefully
acknowledged.}
%%%%%%%%%%%%%%%%

%%%%%%%%%%%%%%
\author{Alexandre Kosyak\corauthref{cor}}
\ead{kosyak01@yahoo.com, kosyak@imath.kiev.ua}
\address{Max-Planck-Institut f\"ur Mathematik, Vivatsgasse 7, D-53111 Bonn, Germany;\\
Institut f\"ur Angewandte Mathematik, Universit\"at Bonn,
Wegelerstr. 6, D-53115 Bonn, Germany; \\
Institute of Mathematics, Ukrainian National Academy of Sciences,
3 Tereshchenkivs'ka, Kyiv, 01601, Ukraine,\\
E-mail: kosyak01@yahoo.com, kosyak@imath.kiev.ua \\
tel.: 38044 2346153 (office), 38044 5656758 (home), fax: 38044
2352010 }

\corauth[cor]{Corresponding author}

\newpage

\begin{abstract}
We give the criteria for the irreducibilty, the Schur irreducibility
and the indecomposability of the set of two $n\times n$ matrices
$\Lambda_n$ and $A_n$ in terms of the subalgebra associated with the
``support" of the matrix $A_n$, where $\Lambda_n$ is a diagonal
matrix with different non zeros eigenvalues and $A_n$ is an
arbitrary one. The list of all maximal subalgebras of the algebra
${\rm Mat}(n,{\mathbb C})$ and the list of the corresponding
invariant subspaces connected with these two matrices is also given.
The properties of the corresponding subalgebras are expressed in
terms of the oriented graphs associated with the support of the
second matrix.

For arbitrary $n$ we describe all minimal subsets of the
elementary matrices $E_{km}$ that generate the algebra ${\rm
Mat}(n,{\mathbb C})$.
\end{abstract}

\end{frontmatter}
\sloppy
\newcommand{\tr}{\mathrm{tr}\,}
\newcommand{\rank}{\mathrm{rank}\,}
\newcommand{\diag}[1]{\mathrm{diag}\,(#1)}
\renewcommand{\Im}{\mathrm{Im}\,}
\maketitle
\newpage
%%%%%%%%%%%
\tableofcontents
%%%%%%%%%%%%%%%%
\section{Introduction} In the {\it representation theory} of different
objects (groups, rings, algebras etc.) the problem of the {\it
ireducibility} of the concrete {\it representations} (modules)
sometimes reduces to the irreducibility of the algebra, generated
by two operators or by two matrices if the representation is
finite dimensional.

In the case of the discrete  group generated by two elements this
is exactly the problem one need to solve. The most popular
examples are the following: the free group ${\mathbb F}_2$
generated by two elements, the Artin braid group $B_3$ on three
strands, the group ${\rm PSL}(2,{\mathbb Z})={\rm SL}(2,{\mathbb
Z})/\pm 1$.

We give the criteria of the {\it irreducibility} and the {\it
Schur irreducibility} (see below the definitions) of the set of
two complex $n\times n$ matrices $\Lambda_n$ and $A_n$ in terms of
the ``support" of the matrix $A_n$, where $\Lambda_n$ is a
diagonal matrix with distinct nonzero eigenvalues and $A_n$ is an
arbitrary one (Theorem \ref{criter123}). The list of {\it all
invariant subspaces } for this two matrices is also given (Theorem
\ref{t.list,S(n)}).

This criterion allows us to study completely in \cite{KosAlb07q}
the irreducibility of some family of representations depending on
the parameters of the braid group $B_3$ in any dimensions.

There are three different notions connected with the {\it
irreducibility} of the representations $T$ of a group $G$ in a
complex vector space $V$
$$
G\ni g\mapsto T_g\in {\rm GL}(V),
$$
where ${\rm GL}(V)$ is the group of  invertible linear operators
on a space $V$. They are the following: 1) {\it irreducible}, 2)
{\it Schur irreducible}, 3) {\it indecomposable}.
\begin{df}  We say that a representation is {\rm
irreducible}  (resp. {\rm  Schur irreducible}) if {\it there are
no nontrivial invariant} closed {\rm subspaces} for all operators
of the representation (resp.  {\rm there are no nontrivial}
bounded {\rm operators commuting} with all operators of the
representation). The representation is {\rm indecomposable} if
{\rm it can not be presented as the direct sum} of the
subrepresentations.
\end{df}
\begin{rem}
\label{dif_irr} It is well known that the relations between the
mentioned notions when the space $V$ is finite dimensional are as
follows: $1)\Rightarrow 2) \Rightarrow 3)$.
\end{rem}
%%
\iffalse
%%
\begin{lem}
\label{dif_irr} We have $1)\Rightarrow 2) \Rightarrow 3)$.
\end{lem}
\begin{pf} We prove that $1)\Rightarrow 2)$. Indeed let we have a
nontrivial operator $A\in B(V)$ such that $[T_g,A]=0\,\, \forall
g\in G$. Let $Sp(A)$ be the spectra of the operator $A$. Since
$Sp(A)\not=\emptyset$ we take some $\lambda\in Sp(A)$ and show
that the subspace $V_\lambda=\{f\in V\mid Af=\lambda f\}$ is
invariant. Indeed let $f\in V_\lambda$ i.e. $Af=\lambda f$. Since
$AT_g=T_gA$ we have $AT_gf=T_gAf=\lambda T_gf$ hence $T_gf\in
V_\lambda$. Since $A$ is nontrivial so $V_\lambda\not=V.$

To prove that $2)\Rightarrow 3)$ we suppose that the
representation is decomposable i.e. holds
$$
V=V_1\oplus V_2\oplus ...\oplus V_k,\quad T_g=T_g^{(1)}\oplus
T_g^{(2)}\oplus ...\oplus T_g^{(k)}.
$$
Let us take $A=\lambda_1 I_1\oplus \lambda_2 I_2\oplus ...\oplus
\lambda_k I_k$ with different $\lambda_i\in {\mathbb C}$.
Obviously, $[A,T_g]=0\,\,\forall g\in G$ and the operator $A$ is
not trivial. \qed\end{pf}
%%
\fi
%
\begin{rem} The notions of  the irreducibility and the Schur irrereducibility
coincide for the unitary representation of an arbitrary group $G$
(hence, for an arbitrary representation of a compact group due to
the "Wayl trick" \cite{Weyl46}).
\end{rem}
{\bf Counterexample 1.} $2)\not\Rightarrow 1)$. Let us consider
the subalgebra of  the algebra ${\rm Mat}(2,{\mathbb C})$
consisting of matrices
$$
%\left(\begin{array}{cc}
\left(\begin{smallmatrix}
a&b\\
0&c
\end{smallmatrix}\right)
%\end{array}\right)
, \,\,a,\,\,b,\,\,c\, \in {\mathbb C}.
$$
This subalgebra is subspace reducible (the subspace in ${\mathbb
C}^2$ generated by the vector $(1,0)$ is invariant) but the
algebra is Schur irreducible.

{\bf Counterexample 2.} $3)\not\Rightarrow 2)$. The classical
example of the Schur reducible but the indecomposable
representation of the additive group of ${\mathbb C}$ is as
follows:
\begin{equation}
\label{-2)+3)} {\mathbb C}\ni z\mapsto \left(\begin{smallmatrix}
1&z\\
0&1
\end{smallmatrix}\right)\in {\rm GL}(2,{\mathbb C}).
\end{equation}
\subsection{Irreducibility criteria}
Let ${\rm Mat}(n,{\mathbb C})$ be the algebra of all complex
matrices over the field of complex numbers ${\mathbb C}$ and let
$\Lambda_n$ (resp. $A_n$) be a diagonal (resp. an arbitrary)
matrix in ${\rm Mat}(n,{\mathbb C})$:
$$
\Lambda_n={\rm
diag}(\lambda_1,\lambda_2,...,\lambda_n),\,\,A_n=(a_{km})_{k,m=1}^n\in{\rm
Mat}(n,{\mathbb C}).
$$
\begin{df} We call {\rm the  support of the matrix} $A=(a_{km})_{k,m=1}^n$ the
subset of indeces $(k,m)$ for which the corresponding entries
$a_{km}$ are nonzero i.e.
\begin{equation}\label{Supp(A)}
{\rm Supp}(A)=\{(k,m)\in \{1,2,...,n\}^2\mid a_{km}\not=0\}.
\end{equation}
\end{df}
\begin{rem}
\label{mat()-s.irr} It is well known that the algebra ${\rm
Mat}(n,{\mathbb C})$ acting on the space ${\mathbb C}^n$ is
irreducible (and hence {\rm Schur irreducible}).
\end{rem}
{\bf Notation.} Denote by $E_{km}$ the matrix units i.e. the
matrix in which the $(k,m)$ entry is $1$ and the other are zero.
%$E_{km}=(a_{ij})$ such that $a_{ij}=\delta_{ki}\delta_{mj}$ where
Obviously $E_{km}E_{pq}=\delta_{mp}E_{kq},$ where $\delta_{mp}$
are the Kronecker symbols.
\begin{thm}
\label{criter123} Let all eigevnalues $\lambda_k$ of $\Lambda_n$
be distinct and nonzero. Then\\
(i) the family of two matrices $(\Lambda_n,A_n)$ is {\rm
irreducible} if and only if the set $\{E_{km}\mid (k,m)\in {\rm
Supp}(A_n)\}$ generates the algebra ${\rm Mat}(n,{\mathbb C})$;\\
(ii) the family  $(\Lambda_n,A_n)$ is {\rm Schur irreducible} if
and only if the set $\{E_{km}\mid (k,m)\in {\rm Supp}(A_n)\bigcup
{\rm Supp}(A^t_n)\}$ generates the algebra ${\rm Mat}(n,{\mathbb
C})$;\\
(iii) the family  $(\Lambda_n,A_n)$ is {\rm indecomposable} if and
only if it is Shur irreducible.
\end{thm}
For a pair of two complex $n\times n$ matrices $A$ and $B$ the
following problems have been studied. When they have 1) a {\it
common eigenvectors}; 2) a {\it common invariant subspace of
dimension} $k,\,\,(2\leq k<n)$?
\\In 1984 Dan Shemesh \cite{Shem84} shows that the criteria for 1)
is: $\bigcap_{k,l=1}^{n-1}{\rm ker}[A^k,B^l]\not=0$.
%%
\iffalse
%%
\begin{thm}
%(The Shemesh criretion).
Let $A,B\in M(n,\mathbb C)$. Define a subspace ${\mathcal N}(A,B)$
of ${\mathbb C}^n$ by formula
\begin{equation}
\label{N(A,B)} {\mathcal N}(A,B) =\bigcap_{k,l=1}^{n-1}{\rm
ker}[A^k,B^l].
\end{equation}
Then the matrices $A$ and $B$ have a common eigenvectors if and
only if
\begin{equation}
\label{N(A,B)not=0} {\mathcal N}(A,B) \not=0.
\end{equation}
\end{thm}
{\bf Problem CIS}. Let $A$ and $B$ be complex $n\times n$
matrices, $k$ a fixed integer, $2\leq k<n$. How can one check
whether $A$ and $B$ have a common invariant subspace of dimension
$k$ by employing a finite number of arithmetic operation?
%%
\fi
%%%
 In \cite{GeoIkra98}, under the additional assumption that at
least one of the matrix $A$ and $B$ has distinct eigenvalues,
were given some sufficient conditions for 2) in terms of $k$th
{\it compound} matrix $C_k(A)$ and $C_k(B)$ of the matrix $A$ and
$B$ (for definition see e.f. \cite{Gan58}, chapt. I,$\S$ 4).
Namely, 2) holds if the matrix $C_k(A)$ and $C_k(B)$ have a common
invariant vector.
% The $k$th{\it compound} matrix $C_k(A)$ of the matrix $A$ is
%a matrix consisting of all minors of order $k$ of the matrix $A$

The {\it advantage of our approach} is that in the case where one
of the matrices is diagonal, we give the {\it criteria for} 2) in
terms of the support of the second matrix. Namely, 2) holds if and
only if subalgebra generated by $\{E_{km}\mid (k,m)\in {\rm
Supp}(A_n)\}$ is contained in  the subalgebra $s_{\bf i}(n)$
defined by (\ref{s_i(n)}) where ${\bf i}=\{i_1,i_2,...,i_k\}$.
%\subseteq \{1,2,...,n\}$.
 The list of {\it
all invariant subspaces } for this two matrices is also given
(Theorem \ref{t.list,S(n)}). In Section 3 we reformulate Theorems
\ref{criter123} and \ref{t.list,S(n)} in terms of the {\it
directed graph} associated with the support of the second matrix.
It allows us to make use of {\it graph theory} (which is well
developed).

%
\iffalse
%
Let $A$ and $B$ be complex $n\times n$ matrices, solution of this
problem under the additional assumption that at least one of the
matrix $A$ and $B$ has distinct eigenvalues.

\begin{thm}[\cite{GeoIkra98}]
%[\cite{GeoIkra98}]
Assume that $A,B\in M(n,\mathbb C)$, $B$ is nonsingular, and $A$
is such that the $k$th {\rm compound} matrix $C_k(A)$ has distinct
eigenvalues for some $k,\,2\leq k<n$. Then if $C_k(A)$ and
$C_k(B)$ have a common eigenvector, $A$ and $B$ have a common
subspace of dimension $k$.
\end{thm}
%%
%%
\begin{df}  {\rm ( see e.f. \cite{Gan58}, chap I,\S 4)}. Let
$A\in M_n(\mathbb C)$ be a given matrix. We consider all possible
minors $M^\alpha_\beta=M^{i_1i_2...i_k}_{j_1j_2...j_k}$ of $A$ of
order $k (1\leq k\leq n).$ The number of this minors is $N^2$,
where $N=\left(\begin{smallmatrix}
n\\
k
\end{smallmatrix}\right)$.
The $k$th {\it compound matrix} of $A$ is $N\times N$ matrix whose
entries are $M^\alpha_\beta$ arranged e.f. lexicographically in
$\alpha$ and $\beta$. This matrix is designated by $C_k(A)$.
\end{df}
%
\fi
%%

\subsection{Irreducibility}
\begin{pf}
$(i)$ The sufficiency  part $``\Leftarrow "$ is obvious due to
Remark \ref{mat()-s.irr}. Indeed let us denote by $\mathfrak{A}_n$
the algebra generated by matrices $\Lambda_n$ and $A_n$. Since
$\lambda_k$ are distinct  and nonzero  we conclude that $E_{kk}\in
\mathfrak{A}_n,\,\,0\leq k\leq n$. Further, since
$$
E_{kk}A_nE_{mm}=a_{km}E_{km},
$$
we conclude that $E_{km}\in \mathfrak{A}_n$ if $a_{km}\not=0$ i.e.
if $(k,m)\in {\rm Supp}(A_n).$

To prove the necessity part $``\Rightarrow "$ for any fixed
$n=1,2,...,$ let us suppose that the set $\{E_{km}\mid (k,m)\in
{\rm Supp}(A_n)\}$ does not generate the whole algebra ${\rm
Mat}(n,{\mathbb C})$, but only some {\it proper subalgebra} $s(n)$
of the following form
\begin{equation}
\label{s(n)1}
s(n)=%\left
\{ x\in {\rm Mat}(n,{\mathbb C})\mid x=\sum_{(k,m)\in S(n)
}x_{km}E_{km}
%\right
\},
\end{equation}
corresponding to some subset of indices $S(n)\subset
\{1,...,n\}^2$.  We can suppose that this subalgebra is {\it
maximal} proper subalgebra of the form (\ref{s(n)1}). Indeed, if
we can find the invariant subspace $V$ for the maximal subalgebra
hence this subspace  would be also invariant one for any of its
subalgebra. By Theorem \ref{t.list,S(n)} the list of the maximal
proper subalgebras $s(n)$ in ${\rm Mat}(n,{\mathbb C})$ of the
form (\ref{s(n)1}) is:
$$
 s_{\bf i}(n)=\{x=(x_{km})_{k,m=1}^n\in{\rm
Mat}(n,{\mathbb C})\mid x_{km}=0,\,k\in\hat{\bf i},\,\,m\in {\bf
i}\},
$$
where ${\bf i}=\{i_1,i_2,...,i_k\}\subseteq
\{1,2,...,n\},\,\,k\leq n$ and $\hat{\bf i}=\{1,2,...,n\}\setminus
{\bf i}.$

{\bf Notation}. For each $n$ let $V_{i_1i_2...i_k}(n):=\langle
e_{i_1},e_{i_2},...e_{i_k}\rangle$ be the linear subspace in
${\mathbb C}^n$ generated by the vectors
$e_{i_1},e_{i_2},...e_{i_k},\,\,1\leq i_1<i_2<...<i_k\leq n$,
where $e_k=(\delta_{rk})_{r=1}^n\in {\mathbb C}^n,\,\,1\leq k\leq
n.$\\
The subspace $V_{\bf i}(n)\!=\!V_{i_1i_2...i_k}(n)$ is an
invariant subspace  for the algebra $s_{\bf i}(n)$.

$(ii)$ To prove the Schur irreducibility we note that the
commutant $(\Lambda_n)':=\{B\in{\rm Mat}(n,{\mathbb C})\mid
[\Lambda_n,B]=0\} $ of the operator $\Lambda_n$ has the following
form:
$$
(\Lambda_n)'=\left(B\in{\rm Mat}(n,{\mathbb C})\mid B={\rm
diag}(b_k)_{k=1}^n\right)
$$
hence, the relation $[A_n,B]=0$ is equivalent to
\begin{equation}
\label{[A,B]=0]}
 a_{km}b_m=b_ka_{km},\quad 1\leq k,m\leq n.
\end{equation}
 We say that we can {\it weakly connect} $k$ and $m$ where
 $ k,m \in\{1,2,\dots, n\}$ if
$a_{km}\not=0$ or $a_{mk}\not=0$, i.e. $(k,m)\in {\rm Supp}(A_n)$
or $(k,m)\in {\rm Supp}(A_n^t)$. In this case $b_k=b_m$. To show
that all $b_k$ coincide (i.e. that $B=bI$) we should be able to
connect step by step all $k$ and $m$ i.e. for any $(k,m)\in
\{1,2,...,n\}^2$ we need to find a sequence
$(k_r,m_r)_{r=1}^l\subseteq {\rm Supp}(A_n)\bigcup {\rm
Supp}(A^t_n),\,\,$ such that
\begin{equation}
\label{gen.E_km} E_{km}=E_{k_1,m_1}E_{k_r,m_r}...E_{k_l,m_l}.
\end{equation}
This proves the sufficiency part of the second part of the
theorem. We say in this case that the set $\{1,2,...,n\}$ is {\it
weakly connected} (see definition (\ref{s.con.dig})).

To prove the necessity  part  let us suppose that the set
$J=\{1,2,...,n\}$ is not  weakly connected i.e. it consists of $l$
weakly  connected components $J_r$ (i.e. $J=\bigcup_{r=1}^lJ_r$).
In this case $b_k=b_m$ for $k,m\in J_r$ and the operator
$B=\oplus_{r=1}^l b_rI_r$ where $I_r=\sum_{k\in J_r}E_{kk}$,
commute with $A_n$, i.e. $[A_n,B]=0$. Hence, the representation is
Schur reducible. Part $(iii)$ is evident. \qed\end{pf}
\section{Maximal proper subalgebras of ${\rm Mat}(n,{\mathbb C})$}
We give the {\it  list} of all subsets of indices
 $S(n)\subset
\{1,2,...,n\}^2$ such that the {\it subspace} $s(n)\in {\rm
Mat}(n,{\mathbb C})$ defined by
\begin{equation}
\label{s(n)}
s(n)=\{ x\in {\rm Mat}(n,{\mathbb C})\mid x=\sum_{(k,m)\in S(n)
}x_{km}E_{km} \}
\end{equation}
is a {\it maximal proper subalgebra} in  ${\rm Mat}(n,{\mathbb
C})$.
\begin{thm}
\label{t.list,S(n)} The list of all maximal proper subalgebras
$s(n)$ in ${\rm Mat}(n,{\mathbb C})$ of the form (\ref{s(n)}) is
as follows
\begin{equation}
\label{s_i(n)}
 s_{\bf i}(n)=\{x=(x_{km})_{k,m=1}^n\in{\rm
Mat}(n,{\mathbb C})\mid x_{km}=0,\,k\in\hat{\bf i},\,\,m\in {\bf
i}\},
\end{equation}
where ${\bf i}=\{i_1,i_2,...,i_k\}\subset \{1,2,...,n\},\,\,k< n$
is a proper subset and $\hat{\bf i}=\{1,2,...,n\}\setminus {\bf
i}.$ The corresponding invariant subspace is $V_{\bf
i}(n):=V_{i_1i_2...i_k}(n)$.
\end{thm}
%To give another proof...
\begin{pf}
For $n=2$ we have just two subsets $S(2)$, namely
$\{(1,1),(1,2),(2,2)\}$ and $\{(1,1),(2,1),(2,2)\}$. Using the
notation (\ref{s_i(n)}) the corresponding  subalgebras are
$$
s_1(2)=\left(\begin{smallmatrix}
*&*\\
0&*\\
\end{smallmatrix}\right),\quad
s_2(2)=\left(\begin{smallmatrix}
*&0\\
*&*\\
\end{smallmatrix}\right).
$$
The mentioned subalgebras $s_1(2)$  and $s_2(2)$ have respectively
the invariant subspaces: $V_1(2)=\langle e_1=(1,0)\rangle$ and
$V_2(2)=\langle e_2=(0,1)\rangle$.

Fort $n=3$ the list of all maximal proper subalgebras is
\begin{equation}
\label{list.s(3)}
 s(3)_{\bf i}:\quad
\left(\begin{smallmatrix}*&*&*\\
0&*&*\\
0&*&*
\end{smallmatrix}\right),\,\,
\left(\begin{smallmatrix}*&0&*\\
*&*&*\\
*&0&*
\end{smallmatrix}\right)
,\,\,
\left(\begin{smallmatrix}*&*&0\\
*&*&0\\
*&*&*
\end{smallmatrix}\right);\,\,
\left(\begin{smallmatrix}*&0&0\\
*&*&*\\
*&*&*
\end{smallmatrix}\right),\,\,
\left(\begin{smallmatrix}*&*&*\\
0&*&0\\
*&*&*
\end{smallmatrix}\right)
,\,\,
\left(\begin{smallmatrix}*&*&*\\
*&*&*\\
0&0&*
\end{smallmatrix}\right).
\end{equation}
The mentioned subalgebras $s(3)$ have respectively the following
invariant subspaces:
$V_1(3),\,\,V_2(3),\,\,V_3(3);\,\,V_{23}(3),\,\,V_{13}(3)$ and
$V_{12}(3)$.

To obtain the list of subalgebras $s(n+1)$ from the list of $s(n)$
we consider two projectors $P_{n,n+1}^{(0)}$ and $P_{n,n+1}^{(1)}$
defined as  follows
$$P_{n,n+1}^{(r)}:{\rm Mat}(n+1,{\mathbb C})\mapsto {\rm
Mat}(n,{\mathbb C}),$$
$$
\sum_{1\leq k,m\leq n+1}x_{km}E_{km}=x\mapsto
P_{n,n+1}^{(r)}(x)=\sum_{r+1\leq k,m\leq n+r}x_{km}E_{km},
$$
$$
\left(\begin{smallmatrix}
*&*&...&*&*\\
&&...&&\\
*&*&...&*&*\\
*&*&...&*&*\\
\end{smallmatrix}\right)
\stackrel{P_{n,n+1}^{(0)}}{\rightarrow}
 \left(\begin{smallmatrix}
*&*&...&*&0\\
&&...&&\\
*&*&...&*&0\\
0&0&...&0&0\\
\end{smallmatrix}\right),\quad
\left(\begin{smallmatrix}
*&*&...&*&*\\
*&*&...&*&*\\
&&...&&\\
*&*&...&*&*\\
*&*&...&*&*\\
\end{smallmatrix}\right)
\stackrel{P_{n,n+1}^{(1)}}{\rightarrow}
 \left(\begin{smallmatrix}
0&0&...&0&0\\
0&*&...&*&*\\
&&...&&\\
0&*&...&*&*\\
\end{smallmatrix}\right).
$$
{\bf Notation}.  For an arbitrary subset of indices   ${\bf
i}=\{i_1,i_2,...,i_k\}\subseteq \{1,2,...,n\}$ let us denote by
$s^{(r)}_{\bf i }(n)=(P_{n,n+1}^{(r)})^{-1}(s_{\bf i }(n))$ the
corresponding subspace in the algebra ${\rm Mat}(n+1,{\mathbb
C})$, where we denote by $A^{-1}(H_0)=\{x\in H_1\mid Ax\in H_0\}$
the preimage of the subset $H_0\subset H_2$ for an operator
$A:H_1\to H_2$.

Let us show how to obtain the list $s(3)$ from the list $s(2)$.
Since the algebra $s(3)$ is contained in the space
$s^{(r)}(2)=(P_{2,3}^{(r)})^{-1}(s(2))$ for $r=0,1$ we get
$$
s_1(2)=\left(\begin{smallmatrix}
*&*\\
0&*\\
\end{smallmatrix}\right)\leftarrow
\left(\begin{smallmatrix}
*&*&\times\\
0&*&\times\\
\times&\times&\times\\
\end{smallmatrix}\right)
,\, \left(\begin{smallmatrix}
\times&\times&\times\\
\times&*&*\\
\times&0&*\\
\end{smallmatrix}\right);\,\,
%%%
s_2(2)=\left(\begin{smallmatrix}
*&0\\
*&*\\
\end{smallmatrix}\right)\leftarrow\left(\begin{smallmatrix}
*&0&\times\\
*&*&\times\\
\times&\times&\times\\
\end{smallmatrix}\right),\,\,
\left(\begin{smallmatrix}
\times&\times&\times\\
\times&*&0\\
\times&*&*\\
\end{smallmatrix}\right).
$$
Since $E_{21}=E_{23}E_{31}$ and $E_{32}=E_{31}E_{12}$ we have just
two subalgebras in $s^{(0)}_1(2)$ and two subalgebras in
$s^{(1)}_1(2)$:
$$
\left(\begin{smallmatrix}
*&*&\times\\
0&*&\times\\
\times&\times&\times\\
\end{smallmatrix}\right)
\to
\left(\begin{smallmatrix}
*&*&*\\
0&*&*\\
0&*&*\\
\end{smallmatrix}\right),\,\,
\left(\begin{smallmatrix}
*&*&*\\
0&*&0\\
*&*&*\\
\end{smallmatrix}\right)
;\quad
 \left(\begin{smallmatrix}
\times&\times&\times\\
\times&*&*\\
\times&0&*\\
\end{smallmatrix}\right)\to
\left(\begin{smallmatrix}
*&0&*\\
*&*&*\\
*&0&*\\
\end{smallmatrix}\right),\,\,\left(\begin{smallmatrix}
*&*&*\\
*&*&*\\
0&0&*\\
\end{smallmatrix}\right)
$$
and  since $E_{12}=E_{13}E_{32}$ and $E_{23}=E_{21}E_{13}$ we have
only two subalgebras in $s^{(0)}_2(2)$ and two subalgebras in
$s^{(1)}_2(2)$:
$$
\left(\begin{smallmatrix}
*&0&\times\\
*&*&\times\\
\times&\times&\times\\
\end{smallmatrix}\right)\to
\left(\begin{smallmatrix}
*&0&*\\
*&*&*\\
*&0&*\\
\end{smallmatrix}\right),\,
\left(\begin{smallmatrix}
*&0&0\\
*&*&*\\
*&*&*\\
\end{smallmatrix}\right);
\,\, \left(\begin{smallmatrix}
\times&\times&\times\\
\times&*&0\\
\times&*&*\\
\end{smallmatrix}\right)\to
\left(\begin{smallmatrix}
*&*&0\\
*&*&0\\
*&*&*\\
\end{smallmatrix}\right),\,
\left(\begin{smallmatrix}
*&*&*\\
0&*&0\\
*&*&*\\
\end{smallmatrix}\right).
$$
Finally we obtain the list (\ref{list.s(3)}) of subalgebras
$s(3)$.
We see that
\begin{align}
\label{(2)to(3)} s_{1}^{(0)}(2)\to s_{1}(3),\,s_{13}(3),\quad
&s_{1}^{(1)}(2)\to
s_{2}(3),\,s_{12}(3),\\
s_{2}^{(0)}(2)\to s_{2}(3),\,s_{23}(3),\quad &s_{2}^{(1)}(2)\to
s_{3}(3),\,s_{13}(3).\nonumber
\end{align}
The list of subalgebra $s(4)$ is as follows:
$$
s_{i}(4):\quad
%1
\left(\begin{smallmatrix}*&*&*&*\\
0&*&*&*\\
0&*&*&*\\
0&*&*&*
\end{smallmatrix}\right),\,\,
%2
\left(\begin{smallmatrix}
*&0&*&*\\
*&*&*&*\\
*&0&*&*\\
*&0&*&*
\end{smallmatrix}\right)
,\,\,
%3
\left(\begin{smallmatrix}
*&*&0&*\\
*&*&0&*\\
*&*&*&* \\
*&*&0&*
\end{smallmatrix}\right),\,\,
\left(\begin{smallmatrix}
*&*&*&0\\
*&*&*&0\\
*&*&*&0 \\
*&*&*&*
\end{smallmatrix}\right),\,\,
$$
$$
s_{i_1i_2i_3}(4):\quad
%1
\left(\begin{smallmatrix}*&0&0&0\\
*&*&*&*\\
*&*&*&*\\
*&*&*&*
\end{smallmatrix}\right),\,\,
%2
\left(\begin{smallmatrix}
*&*&*&*\\
0&*&0&0\\
*&*&*&*\\
*&*&*&*
\end{smallmatrix}\right)
,\,\,
%3
\left(\begin{smallmatrix}
*&*&*&*\\
*&*&*&*\\
0&0&*&0 \\
*&*&*&*
\end{smallmatrix}\right),\,\,
\left(\begin{smallmatrix}
*&*&*&*\\
*&*&*&*\\
*&*&*&* \\
0&0&0&*
\end{smallmatrix}\right),\,\,
$$
$$
s_{i_1i_2}(4):\quad
%12
\left(\begin{smallmatrix}
*&*&*&*\\
*&*&*&*\\
0&0&*&*\\
0&0&*&*
\end{smallmatrix}\right)
%13
,\,\, \left(\begin{smallmatrix}
*&*&*&*\\
0&*&0&*\\
*&*&*&* \\
0&*&0&*
\end{smallmatrix}\right)
%14
\left(\begin{smallmatrix}
*&*&*&*\\
0&*&*&0\\
0&*&*&0 \\
*&*&*&*
\end{smallmatrix}\right),\,\,
%23
\left(\begin{smallmatrix}
*&0&0&*\\
*&*&*&*\\
*&*&*&* \\
*&0&0&*
\end{smallmatrix}\right),\,\,
%24
\left(\begin{smallmatrix}
*&0&*&0\\
*&*&*&*\\
*&0&*&0 \\
*&*&*&*
\end{smallmatrix}\right),\,\,
%34
\left(\begin{smallmatrix}
*&*&0&0\\
*&*&0&0\\
*&*&*&* \\
*&*&*&*
\end{smallmatrix}\right).
$$
The corresponding invariant subspaces are $V_i(4),\,\,1\leq i\leq
4;\,\, V_{i_1i_2i_3}(4),\,\,1\leq i_1<i_2<i_3\leq 4;\,\,$ and
$V_{i_1i_2}(4),\,\,1\leq i_1<i_2\leq 4$.
%%%%%%%%%%%%%ind 3-4

To get $s(4)$ from $s(3)$ we show how this works only for two
subalgebras $s_{1}(3)$ and $s_{13}(3)$
$$
\left(\begin{smallmatrix}
*&*&*\\
0&*&*\\
0&*&*
\end{smallmatrix}\right)\leftarrow
\left(\begin{smallmatrix}
*&*&*&\times\\
0&*&*&\times\\
0&*&*&\times\\
\times&\times&\times&\times\\
\end{smallmatrix}\right),\,
\left(\begin{smallmatrix}
\times&\times&\times&\times\\
\times&*&*&*\\
\times&0&*&*\\
\times&0&*&*
\end{smallmatrix}\right),\quad
\left(\begin{smallmatrix}
*&*&*\\
0&*&0\\
*&*&*
\end{smallmatrix}\right)\leftarrow
\left(\begin{smallmatrix}
*&*&*&\times\\
0&*&0&\times\\
*&*&*&\times\\
\times&\times&\times&\times\\
\end{smallmatrix}\right),\,
\left(\begin{smallmatrix}
\times&\times&\times&\times\\
\times&*&*&*\\
\times&0&*&0\\
\times&*&*&*
\end{smallmatrix}\right).
$$
%%%%%%%%%%%%
Since we have only one possibility to obtain $E_{21}$ and
$E_{31}$, namely $E_{21}=E_{24}E_{41}$ and $E_{31}=E_{34}E_{41}$,
we have only two subalgebras in $s^{(0)}_1(3)$ (case $(a)$). The
other cases are treated similarly. In the case
$(b)\,\,s^{(1)}_1(3)$ we have $E_{32}=E_{31}E_{12}$ and
$E_{42}=E_{41}E_{12}$; in the case $(c)\,\,s^{(0)}_{13}(3)$ we
have $E_{21}=E_{24}E_{41}$ and $E_{23}=E_{24}E_{43}$; in the case
$(d)\,\,s^{(1)}_{13}(3)$ we have $E_{32}=E_{31}E_{12}$ and
$E_{34}=E_{31}E_{14}$. Finally we get
$$
(a)\,\,\left(\begin{smallmatrix}
*&*&*&\times\\
0&*&*&\times\\
0&*&*&\times\\
\times&\times&\times&\times\\
\end{smallmatrix}\right)\to
\left(\begin{smallmatrix}
*&*&*&*\\
0&*&*&*\\
0&*&*&*\\
0&*&*&*\\
\end{smallmatrix}\right),\,
\,\,\left(\begin{smallmatrix}
*&*&*&*\\
0&*&*&0\\
0&*&*&0\\
*&*&*&*\\
\end{smallmatrix}\right);
\quad
(b)\,\, \left(\begin{smallmatrix}
\times&\times&\times&\times\\
\times&*&*&*\\
\times&0&*&*\\
\times&0&*&*
\end{smallmatrix}\right)\to
\left(\begin{smallmatrix}
*&0&*&*\\
*&*&*&*\\
*&0&*&*\\
*&0&*&*
\end{smallmatrix}\right),\,
\left(\begin{smallmatrix}
*&*&*&*\\
*&*&*&*\\
0&0&*&*\\
0&0&*&*
\end{smallmatrix}\right),
$$
$$
(c)\,\,\left(\begin{smallmatrix}
*&*&*&\times\\
0&*&0&\times\\
*&*&*&\times\\
\times&\times&\times&\times\\
\end{smallmatrix}\right)\to
\left(\begin{smallmatrix}
*&*&*&*\\
0&*&0&*\\
*&*&*&*\\
0&*&0&*\\
\end{smallmatrix}\right),\,\,
\left(\begin{smallmatrix}
*&*&*&*\\
0&*&0&0\\
*&*&*&*\\
*&*&*&*\\
\end{smallmatrix}\right);\quad
(d)\,\,\left(\begin{smallmatrix}
\times&\times&\times&\times\\
\times&*&*&*\\
\times&0&*&0\\
\times&*&*&*
\end{smallmatrix}\right)\to
\left(\begin{smallmatrix}
*&0&*&0\\
*&*&*&*\\
*&0&*&0\\
*&*&*&*
\end{smallmatrix}\right),\,
\left(\begin{smallmatrix}
*&*&*&*\\
*&*&*&*\\
0&0&*&0\\
*&*&*&*
\end{smallmatrix}\right).
$$
So we have the following relations using (\ref{(2)to(3)}) and the
latter considerations:
\begin{align*}
s_{1}^{(0)}(2)\to s_{1}(3),\,s_{13}(3),\quad
&s_{1}^{(0)}(3)\to s_{1}(4),\,s_{14}(4),\\
s_{1}^{(1)}(2)\to s_{2}(3),\,s_{12}(3),\quad
&s_{1}^{(1)}(3)\to s_{2}(4),\,s_{12}(4),\\
s_{2}^{(0)}(2)\to s_{2}(3),\,s_{23}(3),\quad
&s_{13}^{(0)}(3)\to s_{13}(4),\,s_{134}(4),\\
s_{2}^{(1)}(2)\to s_{3}(3),\,s_{13}(3),\quad &s_{13}^{(1)}(2)\to
s_{24}(4),\,s_{124}(4).
\end{align*}
To guess the general formula for an arbitrary $n$ we note that
$$
s_{14}^{(0)}(4)\to s_{14}(5),\,s_{145}(5).
$$
The similar considerations explains us how to describe all the
subalgebras \\$s(n+1)$ starting from the subalgebras $s(n)$.
Namely we have
$$
s_{\bf i }(n)\leftarrow s^{(0)}_{\bf i }(n),\,\,s^{(1)}_{\bf i
}(n),
$$
$$
s^{(0)}_{\bf i }(n)\to s_{\bf i }(n+1),\,s_{{\bf i }_0}(n+1),\quad
s^{(1)}_{\bf i }(n)\to s_{\bf i+1}(n+1),\,s_{{\bf i}_1}(n+1),
$$
or
%%%%%%

\begin{equation}
{\unitlength=0,8mm
\begin{picture}(120,20)(10,0)
\put(60,20){\makebox(10,5){\mbox{$s_{\bf i }(n)$}}}
\put(58,20){\line(-4,-1){20}} \put(72,20){\line(4,-1){20}}
\put(20,10){\makebox(10,5){\mbox{$s^{(0)}_{\bf i }(n)$}}}
\put(100,10){\makebox(10,5){\mbox{$s^{(1)}_{\bf i }(n)$}}}
\put(14,10){\line(-2,-1){5}} \put(34,10){\line(2,-1){5}}
\put(94,10){\line(-2,-1){5}} \put(114,10){\line(2,-1){5}}
\put(0,0){\makebox(10,5){\mbox{$s_{\bf i }(n+1)$}}}
\put(40,0){\makebox(10,5){\mbox{$s_{{\bf i }_0}(n+1)$}}}
\put(80,0){\makebox(10,5){\mbox{$s_{\bf i+1}(n+1)$}}}
\put(120,0){\makebox(10,5){\mbox{$s_{{\bf i}_1}(n+1)$}}}
% \put(0,0){\makebox(0,0){\mbox{$\circ$}}}
%\put(1,1){\line(1,0){17}} \put(1,-1){\line(1,0){17}}
%\put(20,0){\makebox(0,0){\mbox{$\circ$}}}
%
%\put(160,0){\makebox(20,20){
\end{picture}
}
\end{equation}
%%%%%%%%%
\iffalse
%%%%%%
\begin{equation}\label{ind(n)(n+1)}
\begin{array}{ccccccccc}
&&&&s_{\bf i }(n)&&&&\\
&&s^{(0)}_{\bf i }(n)&&&&s^{(1)}_{\bf i }(n)&&\\
&s_{\bf i }(n+1)&&s_{{\bf i }_0}(n+1)&
&s_{\bf i+1}(n+1)&&s_{{\bf i}_1}(n+1)&\\
\end{array}
\end{equation}
\fi
%%%%%%
 where for ${\bf i}=\{i_1,i_2,...,i_k\}$ we write ${\bf
i+1}=\{i_1+1,i_2+1,...,i_k+1\},\,\,{\bf i}_0={\bf i}\cup \{n+1\}$
and ${\bf i}_1={\bf i+1}\cup \{1\}$.
%%%%%%%%%%%%%%%%
\qed\end{pf}
\section{Generating sets, maximum subalgebra and the graph theory}

\begin{df}
We say that a subset $G\subset \{1,2,...,n\}^2$ is {\rm generating
subset} if the set of matrices
$$
\{E_{km}\mid (k,m)\in G\}
$$
{\it generates the algebra} ${\rm Mat}(n,{\mathbb C})$.
\end{df}
We would like to describe the {\it minimal} generating subsets $
G$ in terms of the graphs. It would be nice also to find the {\it
complete list} $G(n)$ of the {\it minimal} generating subsets in
$\{1,2,...,n\}^2$. We remind several definitions from the graph
theory.
\begin{df}  We associate with any subset
$G\subset\{1,2,...,n\}^2$ an {\rm directed (oriented) graph
(digraph)} $\Gamma$ on $n$ vertices in the usual way: if $(k,m)\in
G$ we draw the edge (arrow, arc) from the vertex $k$ to the vertex
$m$ on the graph.
\end{df}
\begin{df}{\rm The adjacency matrix} of a graph $\Gamma$ is the $n\times n$
matrix $A_\Gamma$, where n is the number of vertices in the graph.
If there is an edge from some vertex $x$ to some vertex $y$, then
the element $a_{x,y}$ is $1$ (or in general the number of $xy$
edges), otherwise it is $0$.
\end{df}
{\bf Notation}. We shall use the same notation for the subset
$G(n)\subseteq \{1,2,\dots,n\}^2$ and for the corresponding
adjacency matrix $A_{G(n)}$ namely, $G(n)=(g_{km})_{k,m=1}^n$
$$
g_{km}= \left\{\begin{array}{ccc} 1,&{\text if\,\,}&(k,m)\in G(n),\\
0,&{\text if\,\,}&(k,m)\not\in G(n).
\end{array}\right.
$$
For $n=2$ we have only one subset $G(2):=\{(1,2),(2,1)\}$:
$$
G(2):=
\left(\begin{smallmatrix}0&1\\
1&0\\
\end{smallmatrix}\right).
$$
For $n=3$ the list of all subsets $G(3)$ is:
\begin{align}\label{G(3)}
G(3):\quad
G_1(3)=\left(\begin{smallmatrix}0&0&1\\
1&0&0\\
0&1&0
\end{smallmatrix}\right),\,\,
%2
%C^2(3)
G_2(3)=\left(\begin{smallmatrix}0&1&0\\
0&0&1\\
1&0&0
\end{smallmatrix}\right),\,\,
G_3(3)=\left(\begin{smallmatrix}0&1&0\\
1&0&1\\
0&1&0
\end{smallmatrix}\right),\\
G_4(3)=\left(\begin{smallmatrix}0&1&1\\
1&0&0\\
1&0&0
\end{smallmatrix}\right),\,\,
G_5(3)=\left(\begin{smallmatrix}0&0&1\\
0&0&1\\
1&1&0
\end{smallmatrix}\right).\nonumber
\end{align}
There are two distinct notions of connectivity in a digraph.
\begin{df}
\label{s.con.dig} A digraph is {\rm weakly connected} if there is
an undirected path between any pair of vertices, and {\rm strongly
connected} if there is a directed path between every pair of
vertices (\cite{Ski90}, Skiena 1990, p. 173).
\end{df}
\begin{lem}
\label{l.v-tg-r} The subset $G$ is {\rm minimal generating} if and
only if the corresponding graph $\Gamma$  is {\rm minimal strongly
connected}.
\end{lem}
\begin{pf} Use (\ref{gen.E_km}). \qed\end{pf}
Now we can reformulate Theorem \ref{criter123} in terms of the
graph $\Gamma_{A_n}$ associated with the support of $G_{A_n}={\rm
Supp}(A_n)$ of the matrix $A_n$.
\begin{thm}
(i)  The family  $(\Lambda_n,A_n)$ is {\rm irreducible} if and
only
if the graph $\Gamma_{A_n}$ is {\rm strongly connected};\\
(ii)  the family $(\Lambda_n,A_n)$ is {\rm Schur irreducible} if
and only if the graph  graph $\Gamma_{A_n}$ is {\rm
weakly connected}; \\
(iii)  the family $(\Lambda_n,A_n)$ is {\rm indecomposable} if and
only if it is {\rm Schur irreducible}.
\end{thm}
%%%%%%%%%%%%%%%%%
Let us  denote by $A_\Gamma=A_G$ the {\it adjacency matrix} of the
graph $\Gamma$ associated  with the a  $G$. We have the following
correspondence:
\begin{equation}
\label{corr:gr-adj} {\rm set\,\,}G\leftrightarrow {\rm
graph\,\,}\Gamma\leftrightarrow {\rm adjacency \,\,matrix
\,\,}A_G=A_\Gamma.
\end{equation}
\begin{df}
For two subsets $G_1,G_2\subseteq \{1,2,...,n\}^2$ define the {\rm
product} $G_3=G_1\circ G_2$ by
\begin{equation}
\label{G_1.G_2} G_1\circ G_2=\{(k,m)\mid (k,p)\in G_1,(p,m)\in
G_2\,\text{ \,\,for\,\,some\,\,} p\}.
\end{equation}
\end{df}
Let us denote for any subset $G\subset\{1,...,n\}^2$ by $g$ the
corresponding subspace
\begin{equation}\label{g(n)}
g=\{ x\in {\rm Mat}(n,{\mathbb C})\mid x=\sum_{(k,m)\in G
}x_{km}E_{km},\,\, x_{km}\in {\mathbb C}\}.
\end{equation}
Let $g_1$ and $g_2$ be the subspaces corresponding (via
(\ref{g(n)})) to two subsets $G_1$ and $G_2$. We define the
product $g_1g_2$ as follows: $g_1g_2={z=xy\mid x\in g_1,y \in
g_2}$. Obviously, we have
$$
g_1g_2=\{ x\in {\rm Mat}(n,{\mathbb C})\mid x=\sum_{(k,m)\in G_3
}x_{km}E_{km} \},
$$
where $G_3=G_1\circ G_2$.

To define correctly the product of two adjacency matrices
$A_{G_1}$ and $A_{G_2}$ corresponding to the subsets $G_1$ and
$G_2$, we assume  that the entries of the matrix $A_{G_i}$, which
are equal to $0$ and $1$, are in the semiring $R$ defined below.
\begin{df}
\label{ring(R)} Denote by $R$ the {\rm  semiring} consisting of
two elements $0$ and $1$ with operations (see \cite{AhoHopUlm76})
\end{df}
\begin{equation}
\label{ring.R} 0+0=0,\quad 0+1=1,\quad 1+1=1,\quad 0\times
0=0,\quad 0\times 1=0,\quad 1\times 1=1,
\end{equation}
\begin{equation}\label{}
\text{We\quad have}\quad A_{G_1}A_{G_2}= A_{G_1\circ G_2}.
\end{equation}
\begin{lem}
\label{G(n)circ} The set $G$ generates the algebra ${\rm
Mat}(n,{\mathbb C})$ if and only if the powers $G^k=G\circ...\circ
G ,\,\,k=1,2,...,n$ cover the set $\{1,2,...,n\}^2$.
\end{lem}
%%%%%%%%%%%
Using Theorem \ref{t.list,S(n)} and Lemma  \ref{l.v-tg-r} we get
\begin{lem}
(i) The number $\sharp(s(n))$ of the maximal proper subalgebras is
equal to
$$
\sharp(s(n))=\sum_{r=1}^{n-1}C_n^r=2^n-2,
$$ the number of ordered subsets of
the set $\{1,2,...,n\}$ of the length between $1$ and $n-1$;\\
(ii) the  number $\sharp(G(n))$ of the generating subset $G(n)$ is
equal to the number of the minimal strongly connected graphs with
$n$ labeled vertices.
\end{lem}
{\bf Problem 1}. To find the number $\sharp(G(n))$ of the minimal
strongly connected digraphs with $n$ labeled vertices for any $n$
(see appendix).
%%%%%%%%%
\section{Appendix, some examples}
%%%%%
{\bf Notation}. For the sake of shortness we use as before the
same notations for the set $G(n)$ and for the corresponding
adjacency matrix $A_{G(n)}$. We shall denote both by  $G(n)$.

{\bf Example 1}. We show, using Lemma \ref{G(n)circ} , that the
set $G(2)$ and sets $G(3)$ from the list (\ref{G(3)}) are
generating.  We use firstly Lemma \ref{G(n)circ} (recall the
definition \ref{ring(R)}).\\
For $n=2$ the set $G(2)$ is obviously unique. We get
$$
G(2)=\left(\begin{smallmatrix}
0&1\\
1&0
\end{smallmatrix}\right),\quad G^2(2)=\left(\begin{smallmatrix}1&0\\
0&1
\end{smallmatrix}\right) \Rightarrow G(2)\bigcup
G(2)^2=\{1,2\}^2.
$$
For $n=3$ we have $G_1(3)\bigcup G_1^2(3)\bigcup
G_1^2(3)=\{1,2,3\}^2$. Indeed
$$
G_1(3)=\left(\begin{smallmatrix}0&0&1\\
1&0&0\\
0&1&0
\end{smallmatrix}\right),\,\,G_1^2(3)=\left(\begin{smallmatrix}0&1&0\\
0&0&1\\
1&0&0
\end{smallmatrix}\right),\,\,G_1^3(3)=\left(\begin{smallmatrix}1&0&0\\
0&1&0\\
0&0&1
\end{smallmatrix}\right),
$$
$$
G_2(3)=G_1^2(3),\,\,G_2^2(3)=G_1^4(3)=G_1(3)
,\,\,G_2^3(3)=G_1^6(3)=G_1^3(3),
$$
$$
G_3(3)=\left(\begin{smallmatrix}
0&1&0\\
1&0&1\\
0&1&0
\end{smallmatrix}\right),\,\,
G_3^2(3)=\left(\begin{smallmatrix}
1&0&1\\
0&1&0\\
1&0&1
\end{smallmatrix}\right)\Rightarrow
 G_3(3)\bigcup
G_3^2(3)=\{1,2,3\}^2,
$$
$$
G_4(3)=\left(\begin{smallmatrix}0&1&1\\
1&0&0\\
1&0&0
\end{smallmatrix}\right)
,\,\, G_4^2(3)=\left(\begin{smallmatrix}
1&0&0\\
0&1&1\\
0&1&1
\end{smallmatrix}\right)\Rightarrow
 G_4(3)\bigcup
G_4^2(3)=\{1,2,3\}^2,
$$
$$
G_5(3)=\left(\begin{smallmatrix}0&0&1\\
0&0&1\\
1&1&0
\end{smallmatrix}\right),\,\,
G_5^2(3)=\left(\begin{smallmatrix}
1&1&0\\
1&1&0\\
0&0&1
\end{smallmatrix}\right)\Rightarrow
 G_5(3)\bigcup
G_5^2(3)=\{1,2,3\}^2.
$$
If we consider the corresponding directed graphs the proof becomes
evident.
%%%%%%%%%%%%
{\bf Example 2}. The list of the non-isomorphic graphs and the
corresponding adjacency matrices for $n=1,2,3$:\\
%% N1
{\unitlength=0,25mm
\begin{picture}(200,5)(0,0)
\put(0,0){\makebox(0,0){\mbox{$\circ$}}}
\put(170,0){\makebox(0,0){\mbox{$(1)$}}}
\end{picture}
}

%%N2
{\unitlength=0,25mm
\begin{picture}(200,5)(0,0)
\put(0,0){\makebox(0,0){\mbox{$\circ$}}} \put(1,1){\line(1,0){17}}
\put(1,-1){\line(1,0){17}}
\put(20,0){\makebox(0,0){\mbox{$\circ$}}}
\put(160,0){\makebox(20,20){\mbox{$
 \left(\begin{smallmatrix}
0&1\\
1&0
\end{smallmatrix}\right)
$}}}
\end{picture}
}

%%N3
{\unitlength=0,25mm
\begin{picture}(200,40)(20,0)
\put(20,0){\makebox(0,0){\mbox{$\circ$}}}
\put(60,0){\makebox(0,0){\mbox{$\circ$}}}
\put(40,30){\makebox(0,0){\mbox{$\circ$}}}
\put(20,0){\line(2,3){20}} \put(20,0){\line(1,0){40}}
\put(60,0){\line(-2,3){20}}
%%%%%%%%%%%%%%%%%
\put(70,0){\makebox(0,0){\mbox{$\circ$}}}
\put(110,0){\makebox(0,0){\mbox{$\circ$}}}
\put(90,30){\makebox(0,0){\mbox{$\circ$}}}
\put(71,1){\line(1,0){39}} \put(71,1){\line(2,3){19}}
\put(71,-1){\line(1,0){39}} \put(70,3){\line(2,3){19}}
\put(160,5){\makebox(20,20){\mbox{$ \left(\begin{smallmatrix}
0&1&0\\
0&0&1\\
1&0&0
\end{smallmatrix}\right)
$}}}
\put(230,5){\makebox(20,20){\mbox{$ \left(\begin{smallmatrix}
0&1&1\\
1&0&0\\
1&0&0
\end{smallmatrix}\right)
$}}}
\end{picture}
}
%%%%%%%%%%%%%%%%%%%%%N4

{\unitlength=0.4mm
\begin{picture}(200,40)(0,0)
\put(0,0){\makebox(0,0){\mbox{$\circ$}}}
\put(20,0){\makebox(0,0){\mbox{$\circ$}}}
\put(0,20){\makebox(0,0){\mbox{$\circ$}}}
\put(20,20){\makebox(0,0){\mbox{$\circ$}}}
\put(0,0){\line(1,0){20}} \put(0,20){\line(1,0){20}}
\put(0,0){\line(0,1){20}} \put(20,0){\line(0,1){20}}
%%%%%%%
\put(30,0){\makebox(0,0){\mbox{$\circ$}}}
\put(50,0){\makebox(0,0){\mbox{$\circ$}}}
\put(30,20){\makebox(0,0){\mbox{$\circ$}}}
\put(50,20){\makebox(0,0){\mbox{$\circ$}}}
\put(30,0){\line(1,1){20}} \put(30,0){\line(0,1){20}}
\put(30,20){\line(1,0){20}} \put(50,0){\line(0,1){20}}
\put(50,0){\line(-1,1){20}}
%%%%%%%%%%
\put(60,0){\makebox(0,0){\mbox{$\circ$}}}
\put(80,0){\makebox(0,0){\mbox{$\circ$}}}
\put(60,20){\makebox(0,0){\mbox{$\circ$}}}
\put(80,20){\makebox(0,0){\mbox{$\circ$}}}
\put(60,0){\line(1,1){20}} \put(60,0){\line(0,1){20}}
\put(60,20){\line(1,0){19}} \put(79,1.5){\line(0,1){19}}
\put(81,1.5){\line(0,1){19}}
%%%%%%%%%%%%
\put(90,0){\makebox(0,0){\mbox{$\circ$}}}
\put(110,0){\makebox(0,0){\mbox{$\circ$}}}
\put(90,20){\makebox(0,0){\mbox{$\circ$}}}
\put(110,20){\makebox(0,0){\mbox{$\circ$}}}
\put(89,2){\line(0,1){17}} \put(90,2){\line(0,1){17}}
\put(109,2){\line(0,1){17}} \put(110,2){\line(0,1){17}}
\put(90.5,19.5){\line(1,0){17.5}} \put(90.5,21){\line(1,0){18}}
%%%%%
\put(120,0){\makebox(0,0){\mbox{$\circ$}}}
\put(140,0){\makebox(0,0){\mbox{$\circ$}}}
\put(120,20){\makebox(0,0){\mbox{$\circ$}}}
\put(140,20){\makebox(0,0){\mbox{$\circ$}}}
\put(119,2){\line(0,1){17}} \put(120,2){\line(0,1){17}}
\put(120,2){\line(1,1){18}} \put(121,1){\line(1,1){18}}
\put(121,1){\line(1,0){17}}\put(121,0){\line(1,0){17}}
\end{picture}
}
$$
\hskip -4.0cm \left(\begin{smallmatrix}
0&1&0&0\\
0&0&1&0\\
0&0&0&1\\
1&0&0&0
\end{smallmatrix}\right)\,
\left(\begin{smallmatrix}
0&1&0&0\\
0&0&1&0\\
1&0&0&1\\
0&1&0&0
\end{smallmatrix}\right)\,
\left(\begin{smallmatrix}
0&1&0&0\\
0&0&1&0\\
1&0&0&1\\
0&0&1&0
\end{smallmatrix}\right)
\left(\begin{smallmatrix}
0&1&0&0\\
1&0&1&0\\
0&1&0&1\\
0&0&1&0
\end{smallmatrix}\right)\,
\left(\begin{smallmatrix}
0&1&1&1\\
1&0&0&0\\
1&0&0&0\\
1&0&0&0
\end{smallmatrix}\right).
$$
%%%%%%%%%%%%%%%%%%%%%%%%%%%%%%%%%%%%%%%%%%%%%%%%%%%%%%%%%
{\bf Example 3}. The list of $G(n)$ and $s(n)$ for $n\leq 4$ (for
$n=4$ only some $G(n)$).
\begin{align}
\nonumber
 G(2):\quad
\left(\begin{smallmatrix}0&1\\
1&0\\
\end{smallmatrix}\right),\,\,s(2):\quad\left(\begin{smallmatrix}
*&*\\
0&*\\
\end{smallmatrix}\right),\,\,
\left(\begin{smallmatrix}
*&0\\
*&*\\
\end{smallmatrix}\right),\\
\nonumber
G(3):\quad
\left(\begin{smallmatrix}0&0&1\\
1&0&0\\
0&1&0
\end{smallmatrix}\right),\,\,
%2
%C^2(3)
\left(\begin{smallmatrix}0&1&0\\
0&0&1\\
1&0&0
\end{smallmatrix}\right),\,\,
%G_3(3)=
\left(\begin{smallmatrix}0&1&0\\
1&0&1\\
0&1&0
\end{smallmatrix}\right),\,\,
%
%G_4(3)=
\left(\begin{smallmatrix}0&1&1\\
1&0&0\\
1&0&0
\end{smallmatrix}\right),\,\,
%G_5(3)=
\left(\begin{smallmatrix}0&0&1\\
0&0&1\\
1&1&0
\end{smallmatrix}\right),\\
\nonumber
s(3):\quad
\left(\begin{smallmatrix}*&*&*\\
0&*&*\\
0&*&*
\end{smallmatrix}\right),\,\,
\left(\begin{smallmatrix}*&0&*\\
*&*&*\\
*&0&*
\end{smallmatrix}\right)
,\,\,
\left(\begin{smallmatrix}*&*&0\\
*&*&0\\
*&*&*
\end{smallmatrix}\right);\,\,
\left(\begin{smallmatrix}*&0&0\\
*&*&*\\
*&*&*
\end{smallmatrix}\right),\,\,
\left(\begin{smallmatrix}*&*&*\\
0&*&0\\
*&*&*
\end{smallmatrix}\right)
,\,\,
\left(\begin{smallmatrix}*&*&*\\
*&*&*\\
0&0&*
\end{smallmatrix}\right),\\
\nonumber
 G(4):\,\,
%4.1
 \left(\begin{smallmatrix}
0&0&0&1\\
1&0&0&0\\
0&1&0&0\\
0&0&1&0
\end{smallmatrix}\right)
,\,\,
%4.2
\left(\begin{smallmatrix}0&1&0&0\\
0&0&1&0\\
0&0&0&1\\
1&0&0&0
\end{smallmatrix}\right),\,\,
%4.3
 \left(\begin{smallmatrix}
0&1&0&0\\
0&0&0&1\\
1&0&0&0\\
0&0&1&0
\end{smallmatrix}\right)
%4.4
 \left(\begin{smallmatrix}
0&0&1&0\\
1&0&0&0\\
0&0&0&1\\
0&1&0&0\\
\end{smallmatrix}\right),\,\,
%4.5
 \left(\begin{smallmatrix}
0&0&1&0\\
0&0&0&1\\
0&1&0&0\\
1&0&0&0\\
\end{smallmatrix}\right),
%4.6
\left(\begin{smallmatrix}
0&0&0&1\\
0&0&1&0\\
1&0&0&0\\
0&1&0&0\\
\end{smallmatrix}\right),
\end{align}
\begin{align}
\nonumber
%
%4.9
\left(\begin{smallmatrix}
0&1&0&0\\
1&0&1&1&\\
0&1&0&0\\
0&1&0&0
\end{smallmatrix}\right),\,\,
\left(\begin{smallmatrix}
0&0&1&0\\
0&0&1&0\\
1&1&0&1\\
0&0&1&0\\
\end{smallmatrix}\right),\,\,
%4.10
\left(\begin{smallmatrix}
0&1&1&1\\
1&0&0&0\\
1&0&0&0\\
1&0&0&0\\
\end{smallmatrix}\right),\,\,
%4.11
\left(\begin{smallmatrix}
0&0&0&1\\
0&0&0&1\\
0&0&0&1\\
1&1&1&0
\end{smallmatrix}\right),
%4.12,
%4.13
\left(\begin{smallmatrix}
0&1&0&0\\
0&0&1&0\\
1&0&0&1\\
0&1&0&0
\end{smallmatrix}\right),\\
%\end{align}
%
%\begin{align}
%
\nonumber
%
%4.14
\left(\begin{smallmatrix}
0&0&1&0\\
1&0&0&1\\
0&1&0&0\\
0&0&1&0
\end{smallmatrix}\right),\quad
%4.15
\left(\begin{smallmatrix}
0&1&0&1\\
0&0&1&0\\
1&0&0&0\\
0&1&0&0\\
\end{smallmatrix}\right)
,\quad
%4.16
\left(\begin{smallmatrix}
0&1&0&0\\
0&0&1&0\\
1&0&0&1\\
0&0&1&0
\end{smallmatrix}\right),\quad
%4.17
 \left(\begin{smallmatrix}
0&1&0&1\\
0&0&1&0\\
1&0&0&0\\
1&0&0&0
\end{smallmatrix}\right)
,\quad
%4.18
 \left(\begin{smallmatrix}
0&1&0&0\\
0&0&1&1\\
1&0&0&0\\
0&1&0&0\\
\end{smallmatrix}\right),\\
%%%%%%%%%%%%%%%%%%%%%%%
%
\nonumber
%
%4.19
 \left(\begin{smallmatrix}
0&1&0&1\\
1&0&0&0\\
0&0&0&1\\
1&0&1&0\\
\end{smallmatrix}\right),
%4.20
\left(\begin{smallmatrix}
0&1&0&1\\
1&0&1&1&\\
0&1&0&0\\
1&0&0&0
\end{smallmatrix}\right),\,\,
%4.21
\left(\begin{smallmatrix}
0&1&0&0\\
1&0&1&0\\
0&1&0&1\\
0&0&1&0\\
\end{smallmatrix}\right),\,\,
%4.22
\left(\begin{smallmatrix}
0&1&1&1\\
1&0&0&0\\
1&0&0&0\\
1&0&0&0\\
\end{smallmatrix}\right),\,\,
%4.23
\left(\begin{smallmatrix}
0&0&1&1\\
0&0&1&0\\
1&1&0&0\\
1&0&0&0
\end{smallmatrix}\right),
%4.24
\left(\begin{smallmatrix}
0&0&0&1\\
0&0&1&1\\
0&1&0&0\\
1&1&0&0\\
\end{smallmatrix}\right),
\end{align}
$$
s_{i}(4):\quad
%1
\left(\begin{smallmatrix}*&*&*&*\\
0&*&*&*\\
0&*&*&*\\
0&*&*&*
\end{smallmatrix}\right),\,\,
%2
\left(\begin{smallmatrix}
*&0&*&*\\
*&*&*&*\\
*&0&*&*\\
*&0&*&*
\end{smallmatrix}\right)
,\,\,
%3
\left(\begin{smallmatrix}
*&*&0&*\\
*&*&0&*\\
*&*&*&* \\
*&*&0&*
\end{smallmatrix}\right),\,\,
\left(\begin{smallmatrix}
*&*&*&0\\
*&*&*&0\\
*&*&*&0 \\
*&*&*&*
\end{smallmatrix}\right),\,\,
$$
$$
s_{i_1i_2i_3}(4):\quad
%1
\left(\begin{smallmatrix}*&0&0&0\\
*&*&*&*\\
*&*&*&*\\
*&*&*&*
\end{smallmatrix}\right),\,\,
%2
\left(\begin{smallmatrix}
*&*&*&*\\
0&*&0&0\\
*&*&*&*\\
*&*&*&*
\end{smallmatrix}\right)
,\,\,
%3
\left(\begin{smallmatrix}
*&*&*&*\\
*&*&*&*\\
0&0&*&0 \\
*&*&*&*
\end{smallmatrix}\right),\,\,
\left(\begin{smallmatrix}
*&*&*&*\\
*&*&*&*\\
*&*&*&* \\
0&0&0&*
\end{smallmatrix}\right),\,\,
$$
$$
s_{i_1i_2}(4):\quad
%12
\left(\begin{smallmatrix}
*&*&*&*\\
*&*&*&*\\
0&0&*&*\\
0&0&*&*
\end{smallmatrix}\right)
%13
,\,\, \left(\begin{smallmatrix}
*&*&*&*\\
0&*&0&*\\
*&*&*&* \\
0&*&0&*
\end{smallmatrix}\right)
%14
\left(\begin{smallmatrix}
*&*&*&*\\
0&*&*&0\\
0&*&*&0 \\
*&*&*&*
\end{smallmatrix}\right),\,\,
%23
\left(\begin{smallmatrix}
*&0&0&*\\
*&*&*&*\\
*&*&*&* \\
*&0&0&*
\end{smallmatrix}\right),\,\,
%24
\left(\begin{smallmatrix}
*&0&*&0\\
*&*&*&*\\
*&0&*&0 \\
*&*&*&*
\end{smallmatrix}\right),\,\,
%34
\left(\begin{smallmatrix}
*&*&0&0\\
*&*&0&0\\
*&*&*&* \\
*&*&*&*
\end{smallmatrix}\right).
$$
{\bf Example 4}.  The number of the minimal strongly connected
digraphs on $n$ labeled vertices for $n\leq 12$.

We recall that a directed graph is {\rm strongly connected} if
there exists a directed path from any vertex to any other, and
{\rm minimal strongly connected} if removing any edge destroys
this property. Using \cite{www.int-num} we get:

\begin{tabular}{|p{1.7cm}|p{0.15cm}|p{0.15cm}|p{0.2cm}|p{0.5cm}|p{0.9cm}|p{1.1cm}|p{1.2cm}|p{1.6cm}|p{2cm}|}
\hline n&1&2&3&4&5&6&7&8&9\\
\hline A130756&1&1&2&5&15&63&288&1526&8627\\
\hline A130768&1&1&5&58&1069&27816&943669&39757264&2010923289\\
\hline
\end{tabular}

\begin{tabular}{|p{2,5cm}|p{2,7cm}|p{3,1cm}|p{3.5cm}|}
\hline 10&11&12&13\\
\hline  52021& 328432& 2160415& \\
\hline  119153235520& 8118839891161& 627023347399296&
54258093698028037\\
\hline
\end{tabular}

A130756  Number of minimally strongly connected digraphs on n
vertices, up to isomorphism: 1, 1, 2, 5, 15, 63, 288, 1526, 8627,
52021, 328432, 2160415,...

More terms from Vladeta Jovovic (vladeta(AT)Eunet.yu), Jul 13
2007.

A130768  Number of minimally strongly connected digraphs on n
labeled vertices: 1, 5, 58, 1069, 27816, 943669, 39757264,
2010923289, 119153235520, 8118839891161, 627023347399296,
54258093698028037,...

{\bf Example 5}. The list of the maximal subalgebras and the
corresponding graphs $n=2,3$:
$$
s(2):\quad\left(\begin{smallmatrix}
*&*\\
0&*\\
\end{smallmatrix}\right),\,\,
\left(\begin{smallmatrix}
*&0\\
*&*\\
\end{smallmatrix}\right),
$$
$$
s(3):\quad
\left(\begin{smallmatrix}*&*&*\\
0&*&*\\
0&*&*
\end{smallmatrix}\right),\,\,
\left(\begin{smallmatrix}*&0&*\\
*&*&*\\
*&0&*
\end{smallmatrix}\right)
,\,\,
\left(\begin{smallmatrix}*&*&0\\
*&*&0\\
*&*&*
\end{smallmatrix}\right);\,\,
\left(\begin{smallmatrix}*&0&0\\
*&*&*\\
*&*&*
\end{smallmatrix}\right),\,\,
\left(\begin{smallmatrix}*&*&*\\
0&*&0\\
*&*&*
\end{smallmatrix}\right)
,\,\,
\left(\begin{smallmatrix}*&*&*\\
*&*&*\\
0&0&*
\end{smallmatrix}\right)
$$
\includegraphics[scale=0.5]{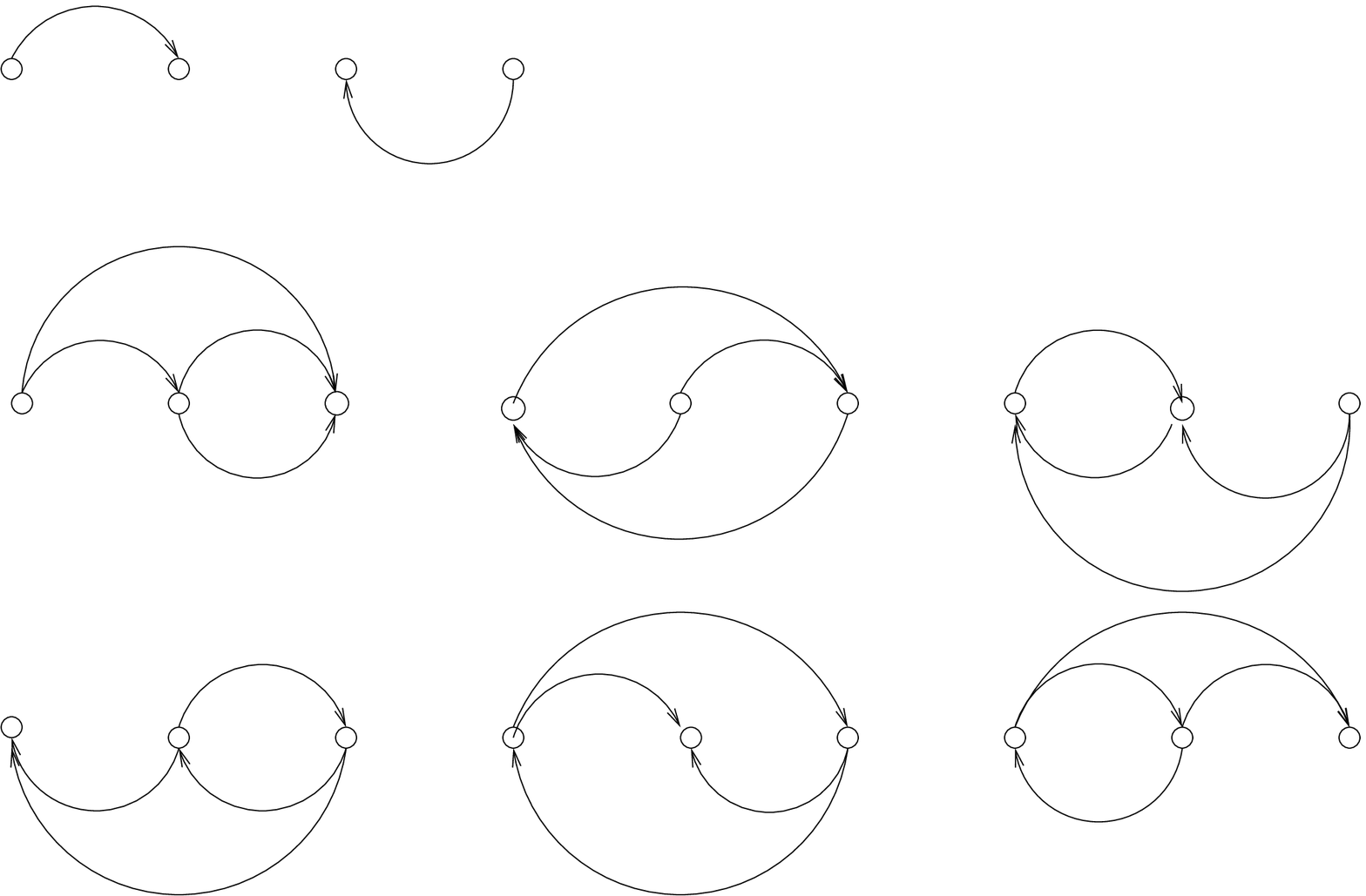}

{\bf Example 6}. The list of the non-isomorphic maximal
subalgebras and the corresponding graphs $n=4$:
$$
%1
s(4)_1=\left(\begin{smallmatrix}*&*&*&*\\
0&*&*&*\\
0&*&*&*\\
0&*&*&*
\end{smallmatrix}\right),\quad
%2
s(4)_{234}=\left(\begin{smallmatrix}*&0&0&0\\
*&*&*&*\\
*&*&*&*\\
*&*&*&*
\end{smallmatrix}\right),\quad
%3
s(4)_{12}=\left(\begin{smallmatrix}
*&*&*&*\\
*&*&*&*\\
0&0&*&*\\
0&0&*&*
\end{smallmatrix}\right),
$$

\includegraphics[scale=0.5]{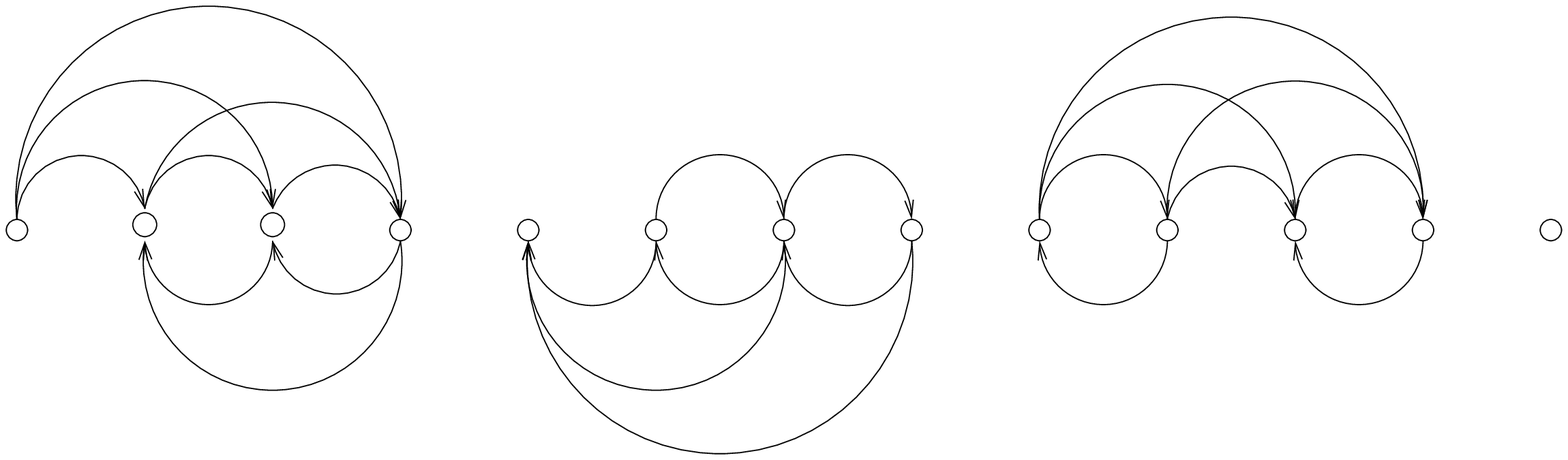}

{\it Acknowledgements.} {The author would like to thank the
Max-Planck-Institute of Mathematics and the Institute of Applied
Mathematics, University of Bonn for the hospitality. The
discussions with Dima Grigoriev, Slava Rabanovich, Volodya
Sergeichuk, Bernard Kt\"otz and Pieter Moree at MPI were very
useful. The author is grateful to Dima Grigoriev for the reference
\cite{AhoHopUlm76} and to Pieter Moree for the reference
\cite{www.int-num} and \cite{Jov07}.  Alexander Schrijver have
kindly explained me
%, due to Pieter Moree,
that the objects I am interested in are exactly the strongly
connected digraphs. The financial support by the DFG project 436
UKR 113/87 is gratefully acknowledged.}

\end{document}